\documentclass[12pt,a4paper]{article}

\usepackage{graphicx}%
\usepackage{amsmath,amsthm,amsfonts,amssymb} 
\usepackage{txfonts,eucal} 
\usepackage{hyperref} 

\newcounter{randnotiz}

\DeclareMathOperator\Char{Char}%
\DeclareMathOperator\Harm{H}%
\DeclareMathOperator\cross{CR}

\newcommand{\cPL}{(\cP,\cL)} 

\renewcommand{\phi}{\varphi}

\newcommand{\cL}{{\mathcal L}}
\newcommand{\cM}{{\mathcal M}}

\newcommand{\cP}{{\mathcal P}}

\newcommand{\cT}{{\mathcal T}}

\newcommand{\bC}{{\mathbb C}}

\newcommand{\bP}{{\mathbb P}}
\newcommand{\bR}{{\mathbb R}}

\newtheorem*{thm*}{Theorem}

\newtheorem*{lem*}{Lemma}
\theoremstyle{definition}

\newtheorem*{defn*}{Definition}
\theoremstyle{remark}

\begin{document}

\author{Hans Havlicek}
\title{Germ\'{a}n Ancochea's work on projectivities, harmonicity preservers and semi-homomorphisms}

\maketitle



\section{Introduction}\label{se:intro}
Our main aim is to analyse three articles of Germ\'{a}n Ancochea (\cite{anco-41a},
\cite{anco-42a}, \cite{anco-47a}) and to describe their impact in algebra and
geometry. Even though most basic mathematical concepts remained unchanged ever
since, terminology, notation, and the style of exposition has undergone
substantial changes. We decided to write our note in today's mathematical
language and to add remarks on the original wording in footnotes. Furthermore,
we provide short expositions of topics that constitute the basis of Ancochea's
research. The surveys \cite{buek+c-95a} and \cite{pick-81a} contain a wealth of
references to the material presented there.

\section{Ancochea's contributions}\label{se:contrib}

\subsection{Projective spaces}\label{subse:proj}
One major theme of Ancochea's contributions is the interplay between synthetic
projective geometry and algebra. Recall that a \emph{projective space} $\cPL$
consists of a set $\cP$ of \emph{points} and a set $\cL$ of subsets of $\cP$
called \emph{lines}. In such a space it is common to introduce the following
terminology. A collection of points is said to be \emph{collinear}, if there
exists a line containing all of them. A \emph{triangle} is a set of three
non-collinear points. A \emph{side} of a triangle is a line that contains two
distinct points of the given triangle. In terms of these notions the axioms of
a projective space $\cPL$ read as follows:
\begin{enumerate}\itemsep0pt
\item Any two distinct points lie on a unique line.

\item If a line meets two distinct sides of a triangle, not at a point of
    the triangle, then it meets the third side also.

\item Any line contains at least three distinct points.
\end{enumerate}
The above definition of a projective space contains no restrictions on its
dimension, which may be finite or infinite.
\par
\begin{figure}[!ht]\unitlength1cm 
  \centering
  \begin{picture}(12.5,3.3)
  \small
    \put(0,0.3){\includegraphics[height=8\unitlength]{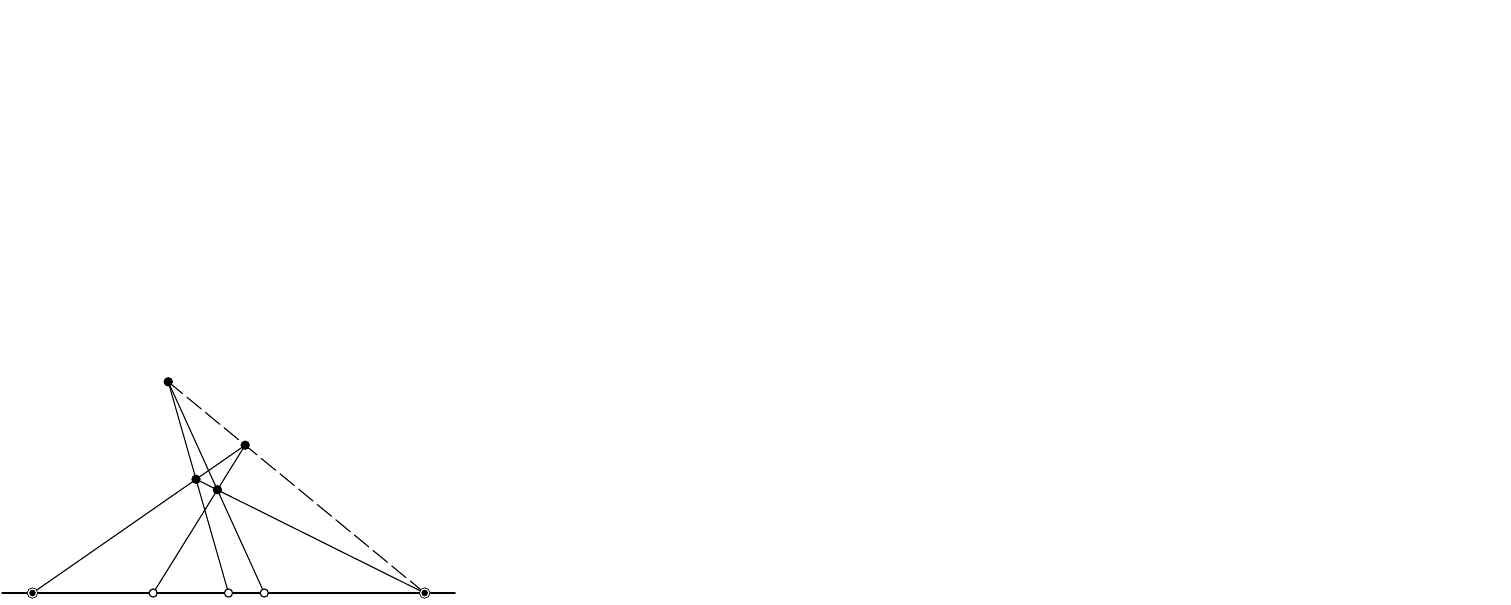}}
    \put(0.0,0.55){$L$}
    \put(0.325,0.0){$0$}
    \put(5.475,0.0){$\infty$}
    \put(1.95,0.0){$y$}
    \put(2.9,0.0){$x$}
    \put(3.4,0.0){$x+y$}
    \put(6.5,0.3){\includegraphics[height=8\unitlength]{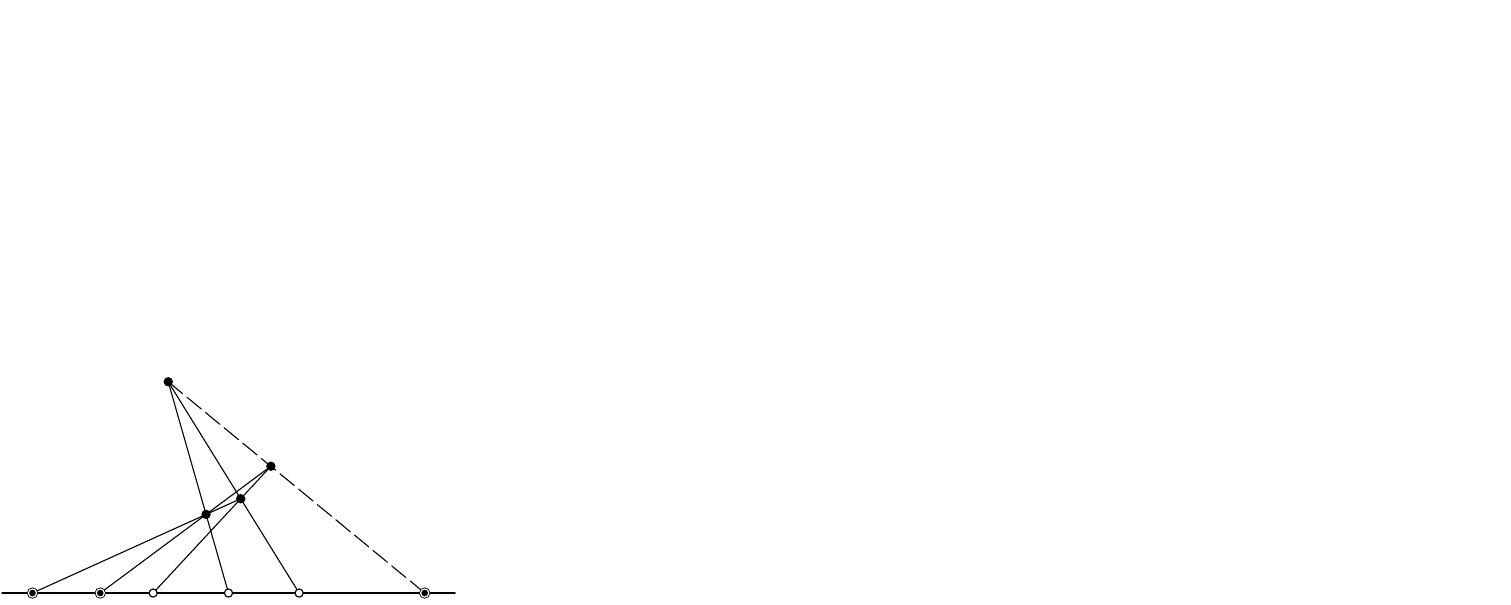}}
    \put(6.5,0.55){$L$}
    \put(6.825,0.0){$0$}
    \put(7.725,0.0){$1$}
    \put(11.975,0.0){$\infty$}
    \put(8.45,0.0){$y$}
    \put(9.4,0.0){$x$}
    \put(10.375,0.0){$x\cdot y$}
  \end{picture}
    \caption{Addition and multiplication}\label{fig:add+mult}
\end{figure}

We now assume, until further notice, $\cPL$ to be a projective space of
dimension $>2$, that is, it contains at least two disjoint lines. Such a space
has a remarkable property. Upon choosing a line $L$ and three distinct
\emph{points of reference} on $L$, which are labelled as $0$, $1$, $\infty$,
the set $K:=L\setminus\{\infty\}$ can be made into a field\footnote{As is
customary among geometers, multiplication in a field is not assumed to be
commutative.} by defining the sum $x+y$ and the product $x\cdot y$ of $x,y\in
K$ in a purely geometric manner. Starting with $x$ and $y$ one has to draw
several auxiliary points and lines in some (projective) plane\footnote{A
\emph{projective plane} is a two-dimensional projective space. When speaking of
subspaces of a projective space, we usually drop the adjective ``projective''.}
containing $L$ according to Fig.~\ref{fig:add+mult}.\footnote{The classical
example of a three-dimensional projective space arises from the
three-dimensional Euclidean space by adding ``points and lines at infinity'' in
an appropriate way. If the dashed line in Fig.~\ref{fig:add+mult} coincides
with a line at infinity and $0$ is not at infinity, then there is an elementary
interpretation. The sum $x+y$ is the image of $x$ under the translation taking
$0$ to $y$. The product $x\cdot y$ is the image of $x$ under the homothety with
centre $0$ taking $1$ to $y$.}
The distinguished points $0$ and $1$ are the neutral elements with respect to
addition and multiplication in the field $K$, which will be called an
\emph{underlying field}\label{page:underlying} of $\cPL$. All underlying fields
of a given projective space $\cPL$ of dimension $>2$ are isomorphic.
\par
\begin{figure}[!ht]\unitlength1cm 
\begin{minipage}[t]{0.5\textwidth}
\centering
  \begin{picture}(6,4)
  \small
    \put(0,0.0){\includegraphics[height=8\unitlength]{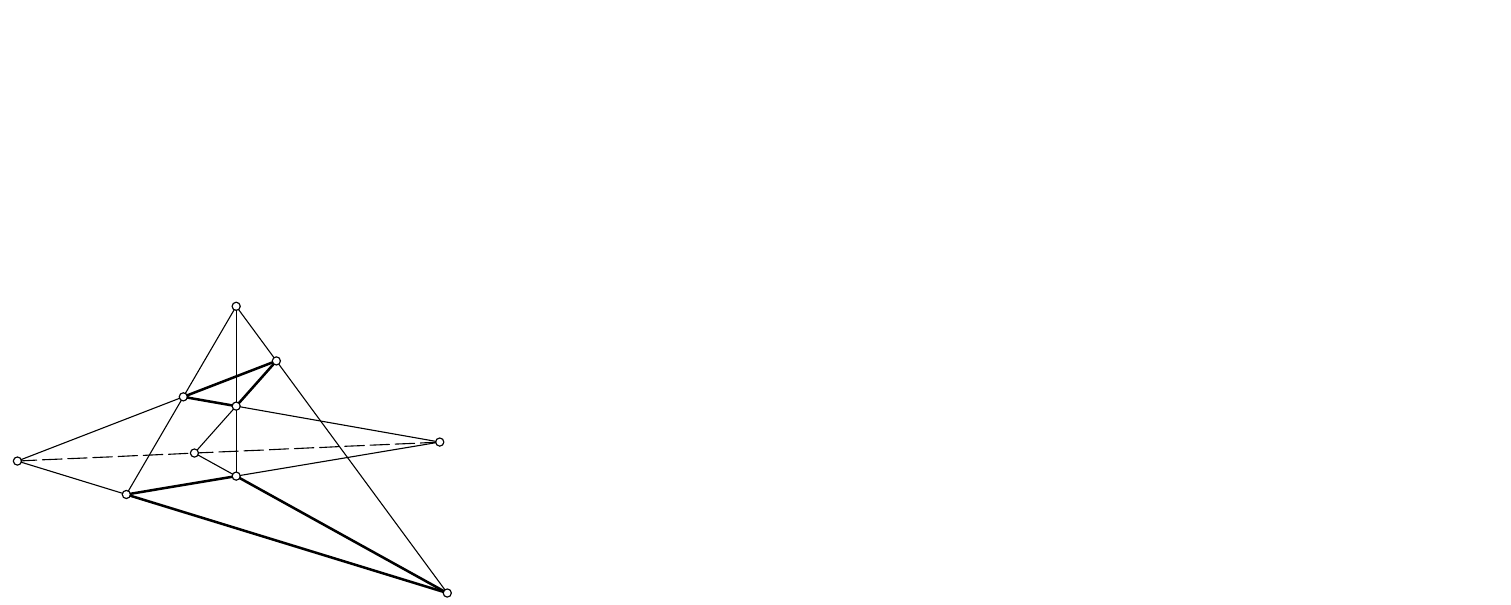}}
    \put(2.8,3.9){$p$}
    \put(2.15,2.9){$a_1$}
    \put(3.175,2.35){$a_2$}
    \put(3.75,3.3){$a_3$}
    \put(1.5,1.1){$b_1$}
    \put(3.05,1.3){$b_2$}
    \put(6.05,0.1){$b_3$}
  \end{picture}
    \caption{Theorem of Desargues}\label{fig:desargues}
\end{minipage}%
\begin{minipage}[t]{0.5\textwidth}
\centering
  \begin{picture}(6,3)
  \small
    \put(0,0.0){\includegraphics[height=8\unitlength]{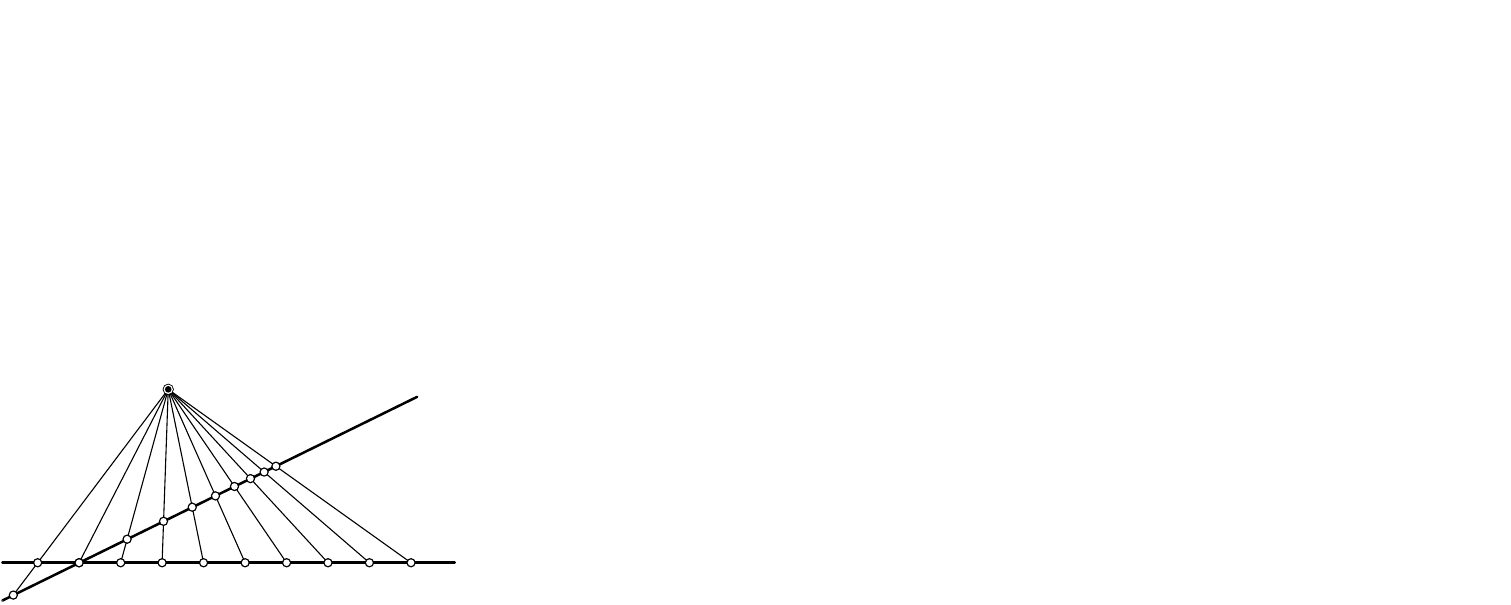}}
    \put(1.9,2.85){$p$}
    \put(5.7,0.70){$L_1$}
    \put(5.1,2.25){$L_2$}
    \put(1.85,0.3){$x$}
    \put(1.8,1.1){$x^\pi$}
  \end{picture}
    \caption{Perspectivity}\label{fig:persp}
\end{minipage}%
\end{figure}

One crucial point is, of course, that the definition of the sum and the product
of points does not depend on the choice of the auxiliary elements. This is due
to the fact that $\cPL$ is a \emph{desarguesian} projective space, that is, in
any of its planes the following \emph{Theorem of Desargues}
(Fig.~\ref{fig:desargues}) holds. \emph{If two triangles $a_1,a_2,a_3$ and
$b_1,b_2,b_3$ of a projective plane are in perspective from a point $p$, then
the intersections of their corresponding sides are collinear.}
\par
Now let $\cPL$ be a projective plane. Then the construction of an underlying
field can be carried out like before provided that $\cPL$ is desarguesian.
Otherwise one arrives at the problem of coordinatisation of a non-desarguesian
projective plane \cite{hugh+p-73a}, \cite{pick-55a}, \cite{pick-75a}, which is
beyond our scope.
\par
Given a projective space $\cPL$ of dimension $\geq 2$ a mapping $\pi\colon
L_1\to L_2\colon x\mapsto x^\pi$ of a line $L_1$ to a line $L_2$ is called a
\emph{perspectivity}, if there exists a point $p\notin L_1\cup L_2$ (called the
\emph{centre} of $\pi$) such that $p,x,x^\pi$ are collinear for all $x\in L$
(Fig.~\ref{fig:persp}).
Any product of finitely many perspectivities is called a
\emph{projectivity}\footnote{This definition goes back to J.-V.~Poncelet.
Ancochea used the name \emph{projectivity in the sense of Poncelet} for such a
mapping.}. Projectivities are bijective. If $(p_1,p_2,p_3)$ and $(q_1,q_2,q_3)$
are triples of distinct points on lines $L$ and $M$, respectively, then there
is at least one projectivity $\psi$ sending $p_1\mapsto q_1$, $p_2\mapsto q_2$,
$p_3\mapsto q_3$. In order to obtain \emph{all} projectivities with this
property it is therefore enough to look for all projectivities of $L$ onto
itself such that $p_1,p_2,p_3$ remain fixed.
\par
In the classical case, where $\cPL$ is three-dimensional and any underlying
field is isomorphic to the field of real numbers, the name \emph{Fundamental
Theorem of Projective Geometry}\label{page:fundament} has been given to the
following result. \emph{The identity is the only projectivity of a line onto
itself with three distinct fixed points.} If $\cPL$ is a desarguesian
projective space of dimension $\geq 2$, then this fundamental theorem remains
valid if, and only if, one underlying field of $\cPL$ is commutative.

\subsection{Sobre el teorema fundamental de la geometria
proyectiva}\label{subse:fundamental}
In this article (\cite{anco-41a}) from 1941, Ancochea considered a projective
space $\cPL$ of dimension $>2$ and referred to the books of A.~N.~Whitehead
\cite[p.~15]{whit-06a} and L.~Bieberbach \cite[p.~5]{bieb-33a} for the notion
of a projective space.\footnote{In both books there are also extra axioms that
force the dimension of a projective space to be three. The formalism used in
\cite{bieb-33a} is different from ours.} Also, he quoted
\cite[Kap.~1,~\S~4]{bieb-33a} for the construction of a field $K$ from a line
of $\cPL$ and recalled the definition and some properties of a projectivity.
\par
After these preliminaries, the main question of the article was posed.
\emph{What can be said about a projectivity of a line $L$ onto itself that
fixes three distinct points?} Ancochea immediately derived a partial answer
from a note by O.~Veblen \cite{vebl-07a}. \emph{If these three points are
labelled as $0,1,\infty$ and are used to construct a field
$K:=L\setminus\{\infty\}$, then such a projectivity restricts to an
automorphism of $K$.} (Let us add that the reasoning in \cite{vebl-07a} relies
on two observations. First, any projectivity can be extended to a
\emph{collineation} of $\cPL$, i.~e. a bijection $\cP\to\cP$ that takes any
line onto a line. Second, any collineation of $\cPL$ fixing the points $0$,
$1$, $\infty$ induces an automorphism of $K$.) Next, he stated the following
complete answer to his question.

\begin{thm*}[{\cite[Teorema fundamental]{anco-41a}}]
Let $\psi$ be a projectivity of a line $L$ with three distinct fixed points. If
these points are chosen as the points of reference $(0,1,\infty)$, then the
image of any point $x\in L$ other than $\infty$ arises by applying an inner
automorphism of $K=L\setminus\{\infty\}$. More precisely, there is an $a\in
K\setminus\{0\}$, which is fixed under $\psi$, such that $x^\psi = a^{-1}xa$
for all $x\in K$. Conversely, any inner automorphism of $K$ extends to a
projectivity of $L$ that fixes the points of reference.
\end{thm*}
Let us sketch Ancochea's elegant proof. One may assume without loss of
generality that $\psi$ is given as a product $\pi_1 \pi_2\cdots\pi_{n+1}$ of
perspectivities and, moreover, that $\pi_{n+1}$ takes the form
$\pi_{n+1}\colon\overline L\to L$ with $0$ being the only common point of $L$
and $\overline L$. By a theorem of F.~Schur \cite{schur-98a}, the projectivity
$\overline\psi:=\pi_1\pi_2\cdots\pi_{n}\colon L\to \overline L$ can also be
written as a product of \emph{at most two} perspectivities. If $\overline \psi$
is a single perspectivity, then $0$, $1$, $\infty$ being fixed yields that
$\psi=\overline\psi \pi_{n+1}$ is the identity. Otherwise, $\psi$ is a product
of three perspectivities. The rest of the proof is contained in
Fig.~\ref{fig:anco}, which is a reproduction\footnote{The author is grateful to
the editors of ``Revista Matem\'{a}tica Iberoamericana'' for granting permission to
include a copy of Ancochea's drawing.} from \cite{anco-41a}.
\begin{figure}[!ht]\unitlength0.80cm 
  \centering
  \begin{picture}(8.640,5.751)
  \small
    \put(0.0,2.5){\begin{minipage}{8.640\unitlength}
                    \centerline{To be included in the final version only.}
                  \end{minipage}
                 }
  \end{picture}
    \caption{Copy of Ancochea's drawing}\label{fig:anco}
\end{figure}
A carefully explained way of how to read this figure\footnote{The points of
reference are labelled as $A_0$, $A_1$, $A_\infty$. In order to illustrate the
case of a non-commutative field $K$, Ancochea decided to ``bend'' one line so
that the point $x$ and its image under $\psi$ (written as $x'$) turn out
different.} in terms of the field $K$ yields that $x^\psi= axa^{-1}$ for an
appropriate $a\in K$. By reversing these arguments the converse can be
established.
\par
Ancochea closed with the observation that the ``classical'' Fundamental Theorem
of Projective Geometry (see page~\pageref{page:fundament}) appears as an
immediate corollary of his more general version. Indeed, the projectivity
$\psi$ reduces to the identity if, and only if, $a$ belongs to the centre of
$K$.
\par
In a subsequent article \cite{anco-42a}, Ancochea expressed his thanks to
F.~Bachmann for pointing out that his \emph{Teorema fundamental} has been
obtained already in 1930, using different methods, by K.~Reidemeister
\cite[p.~137, Satz~4]{reid-30a} (reprint \cite{reid-68a}).

\subsection{Harmonic quadruples}\label{subse:harmonic}
Let $\cPL$ be a projective space of dimension $\geq 2$. In $\cPL$, a
\emph{quadrangle} means a set of four points of a fixed plane no three of which
are collinear. A quadruple $(p_1,p_2,p_3,p_4)$ of points on a line is said to
be \emph{harmonic}, in symbols $\Harm(p_1,p_2,p_3,p_4)$, if there exists a
quadrangle, say $q_1,q_2,q_3,q_4$ with $p_1=q_1q_2\cap q_3q_4$, $p_2=q_2q_3\cap
q_4q_1$, $p_3\in q_1q_3$, and $p_4\in q_2q_4$. (Here we denote the unique line
joining $p$ and $q$ by $pq$ etc.) The first three points $p_1,p_2,p_3$ of a
harmonic quadruple are distinct. Likewise, the fourth point $p_4$ cannot
coincide with $p_1$ or $p_2$. In general, nothing can be said about $p_3$ and
$p_4$ being different or not. However, if $\cPL$ is desarguesian, then harmonic
quadruples are well understood. First of all, the theorem about the
\emph{uniqueness of the fourth harmonic point}\label{page:uniqueness} holds.
\emph{Given three distinct points $p_1,p_2,p_3$ on a line $L$ there is a unique
point $p_4\in L$ such that $\Harm(p_1,p_2,p_3,p_4)$}. Furthermore, an
underlying field $K$ has characteristic two if, and only if, $p_3=p_4$.
\begin{figure}[!ht]\unitlength1cm 
  \centering
  \begin{picture}(6,2.5)
  \small
    \put(0,0.0){\includegraphics[height=8\unitlength]{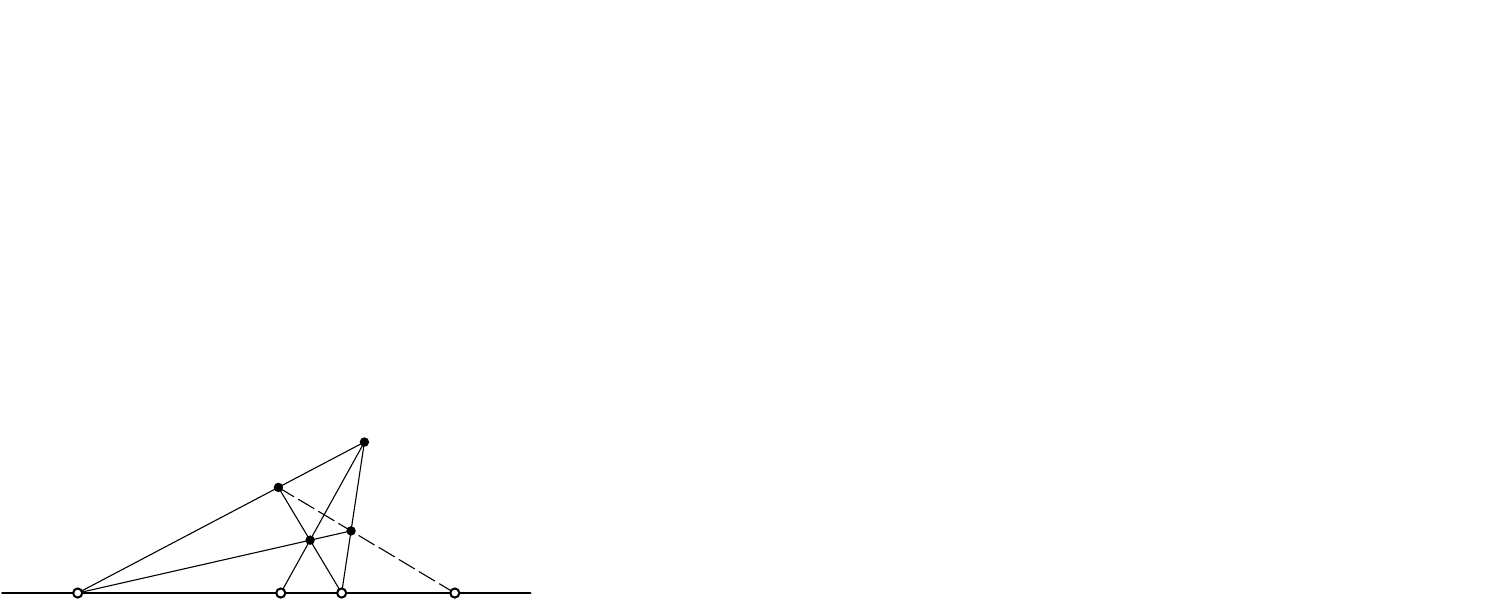}}
    \put(4.7,2.35){$q_1$}
    \put(3.4,1.75){$q_2$}
    \put(3.65,0.9){$q_3$}
    \put(4.75,1.1){$q_4$}
    \put(0.95,-0.2){$p_1$}
    \put(4.35,-0.2){$p_2$}
    \put(3.55,-0.2){$p_3$}
    \put(5.95,-0.2){$p_4$}
  \end{picture}
    \caption{Harmonic quadruple}\label{fig:harm}
\end{figure}
\par
Let us now switch from synthetic to analytic projective geometry. Given any
left vector space $V$ over a field $K$ the \emph{projective space on $V$},
which will be written as $\bP(V)$, has $\bigl\{Kv \mid v\in
V\setminus\{0\}\bigr\}$, i.~e. the set of one-dimensional subspaces of $V$, as
its set of points. A line of $\bP(V)$ is defined as the set of all points that
are contained in some fixed subspace of $V$ with (vector) dimension two. We
shall also use the name \emph{projective space over $K$} for $\bP(V)$ in order
to emphasise the ground field. A projective space $\bP(V)$ has dimension $\dim
V-1$, and it is desarguesian.\footnote{The projective spaces $\bP(V)$ with
$\dim V\geq 3$ are, up to isomorphism, precisely the desarguesian projective
spaces of dimension $\geq 2$.}
\par
If $(p_1,p_2,p_3)$ is a triple of distinct points on a line $L$ of $\bP(V)$,
then there exists at least one pair of linearly independent vectors $(w_1,
w_2)$ such that $p_1=Kw_1$, $p_2=Kw_2$, $p_3=K(w_1+w_2)$. The mapping
\begin{equation}\label{eq:identif}
    L\setminus\{p_1\}\ \to K\colon K(x_1w_1+x_2w_2) \mapsto x_2^{-1}x_1
    \mbox{~~with~} x_1,x_2\in K
    \mbox{~and~} x_2\neq 0
\end{equation}
is bijective. If $K$ is non-commutative, then this bijection depends not only
on the given triple of points but also on the choice of $(w_1,w_2)$. Indeed,
all pairs of vectors with the required property comprise the set $\big\{
(cw_1,cw_2)\mid c\in K\setminus\{0\} \big\}$. When fixing any point
$p_4=K(x_1w_1+x_2w_2)\in L\setminus\{p_1\}$ and calculating the corresponding
element of $K$ with respect to all these pairs in analogy to
\eqref{eq:identif}, one obtains the set
\begin{equation}\label{eq:cross}
    \cross(p_1,p_2,p_3,p_4) :=
    \left\{c(x_2^{-1}x_1)c^{-1}\mid c\in K\setminus\{0\} \right\},
\end{equation}
which is called the \emph{cross ratio} of the quadruple $(p_1,p_2,p_3,p_4)$.
This cross ratio is a \emph{conjugacy class} in $K$ and it depends only on the
given points. Precisely when $x_2^{-1}x_1$ is in the centre of $K$, the above
cross ratio may be considered as the single element $x_2^{-1}x_1\in K$ rather
than a subset of $K$.\footnote{Ancochea never mentions ``cross ratios'' in
\cite{anco-41a}, \cite{anco-42a}, and \cite{anco-47a}. When he wrote these
articles, the notion of ``cross ratio'' has been established for commutative
fields only. The definition in \eqref{eq:cross} was given in 1952 by R.~Baer
\cite[pp.~71--72]{baer-52a} and, in a slightly different setting, in 1948--1949
by E.~Sperner \cite[p.~149]{sper-48a} and \cite[p.~425]{sper-49a}.}
\par
If $\dim V > 2 $, then harmonic quadruples can be characterised as follows:
\begin{equation}\label{eq:cr=-1}
    \cross(p_1,p_2,p_3,p_4) = -1 \;\;\Leftrightarrow \;\;\Harm(p_1,p_2,p_3,p_4) .
\end{equation}
Also, with $(0,1,\infty):=(p_2,p_3,p_1)$ the bijection in \eqref{eq:identif}
turns into an isomorphism of the underlying field $L\setminus\{\infty\}$ onto
$K$. So, all underlying fields of $\bP(V)$ are isomorphic to $K$. If $\dim V=2$
or, in other words, if $\bP(V)$ is a projective line, then \eqref{eq:cr=-1} can
be used as definition of harmonic quadruples. An abstract projective line,
however, has in general not enough ``intrinsic structure'' to define harmonic
quadruples.
\par
The seminal book \emph{Geometrie der Lage} by K.~G.~C.~von Staudt contains a
most remarkable result, namely a characterisation of projectivities in terms of
harmonic quadruples \cite[pp.~49ff.]{stau-47a}. In its original version, which
is confined to projective spaces of dimension three over the real numbers
$\bR$, it reads as follows:

\begin{thm*}[{von Staudt's Theorem}]
Let lines $L_1$ and $L_2$ be given. Then a bijective mapping $\lambda\colon
L_1\to L_2$ is a projectivity precisely when it takes any harmonic quadruple on
$L_1$ to a harmonic quadruple on $L_2$.
\end{thm*}
For this result to be true, two properties of the field $\bR$ are essential.
First, the characteristic $\Char \bR$ is different from two and, second, the
field $\bR$ admits no automorphism other than the identity. This was pointed
out in a note by M.~G.~Darboux \cite{darb-80a} from 1870, where also a small
gap in von Staudt's original proof was closed.
\par
In order to generalise von Staudt's Theorem we adopt the following terminology.
Given lines $L_1$ and $L_2$ in a desarguesian projective space $\cPL$ of
dimension $\geq 2$ a mapping $\lambda\colon L_1\to L_2$ is called a
\emph{harmonicity preserver}\label{page:h-preserver}\footnote{Some authors use
the name \emph{harmonic mapping} instead. We refrain from doing so in order to
avoid confusion with the ``harmonic mappings'' known from differential
geometry.} if it takes any harmonic quadruple on $L_1$ to a harmonic quadruple
on $L_2$.
\par
Harmonicity preservers do not deserve interest if an underlying field of $\cPL$
is of characteristic two, since then a quadruple $(p_1,p_2,p_3,p_4)$ of points
is harmonic if, and only if, $p_1,p_2,p_3$ are distinct collinear points and
$p_3=p_4$. So, in case of characteristic two, a mapping $\lambda\colon L_1\to
L_2$ is a harmonicity preserver precisely when it is injective.
\par
Suppose now that $\cPL$ has an underlying field of characteristic $\neq 2$.
Projectivities are obvious examples of \emph{bijective} harmonicity preservers.
In order to determine all bijective harmonicity preservers between lines $L_1$
and $L_2$ of $\cPL$, it is therefore enough to consider all mappings $\lambda$
of this kind that fix three arbitrarily chosen points on some line $L$. When
labelling the three fixed points as $0$, $1$, $\infty$, the points of
$K:=L\setminus\{\infty\}$ provide an underlying field of $\cPL$, and $\lambda$
restricts to a bijection $\sigma$ of $K$; see page \pageref{page:underlying}.
O.~Schreier and E.~Sperner \cite[pp.~191--194]{schrei+s-35a} used these ideas
and extended the findings of von Staudt and Darboux in the year 1935. \emph{If
the field $K$ is commutative, $\Char K\neq 2$, then any bijective harmonicity
preserver $\lambda\colon L\to L$ fixing the points $0$, $1$, $\infty$ restricts
to an automorphism of $K$ and, conversely, any automorphism of $K$ extends to a
unique bijective harmonicity preserver of the line $L$.}

\subsection{Le th\'{e}or\`{e}me de von Staudt en g\'{e}om\'{e}trie projective
quaternionienne}\label{subse:staudt}
Ancochea started his article \cite{anco-42a} by quoting the result of
O.~Schreier and E.~Sperner that we encountered above. He observed that an
extension thereof to the case of a non-commutative field $K$ seemed to be
unknown. In order to state such an extension he introduced an entirely new
algebraic concept.

\begin{defn*}[{\cite[p.~193]{anco-42a}}] A bijective mapping $\sigma$ of a field $K$ onto itself is
called a \emph{semi-automorphism} if for all $x,y\in K$ the following
conditions hold:
    \begin{align}
  (x+y)^\sigma             &= x^\sigma + y^\sigma ,                  \label{eq:I}\\
  (xy)^\sigma + (yx)^\sigma &= x^\sigma   y^\sigma + y^\sigma x^\sigma .\label{eq:II}
    \end{align}
\end{defn*}
This allowed him to state the following theorem.

\begin{thm*}[{\cite[p.~193]{anco-42a}}]
In a projective space over a field of characteristic $\neq 2$, any bijective
harmonicity preserver\footnote{Ancochea used instead the name
\emph{projectivity in the sense of von Staudt}.} of a line $L$ onto itself with
three distinct fixed points $0$, $1$, $\infty$ restricts to a semi-automorphism
of the field $L\setminus\{\infty\}$.
\end{thm*}
In order to establish the theorem, Ancochea considered the projective space
$\bP(V)$ on some left vector space $V$ over a field $K$ of characteristic $\neq
2$. The given line $L$ therefore corresponds to a two-dimensional subspace $W$
of $V$. There exists a basis $(w_1,w_2)$ of $W$ such that the distinguished
points $0$, $1$, $\infty$ are $Kw_2$, $K(w_1+w_2)$, $Kw_1$,
respectively.\footnote{Ancochea actually used the left vector space $K^2$ and
its canonical basis instead of our $W$ and $(w_1,w_2)$.} The isomorphism
\eqref{eq:identif} allowed him to \emph{identify} the points of the underlying
field $L\setminus\{\infty\}$ with the elements of the field $K$. Next, he
addressed the problem of calculating the fourth harmonic point, say $p$ of
three distinct points $a_1,a_2,a_3\in K$. Based on calculations, which in most
cases were omitted by Ancochea, he arrived at the following results: $p=\infty$
is equivalent to $ a_3=\frac12(a_1+a_2)$. If $p\neq\infty$, letting $p=:a_4\in
K$ gives
\begin{equation}\label{eq:mult}
  a_4 = \left( (a_1-a_3)^{-1} + (a_2-a_3)^{-1} \right)^{-1}
  \left( (a_1-a_3)^{-1}a_1 + (a_2-a_3)^{-1}a_2 \right) .
\end{equation}
Equation \eqref{eq:mult} takes the form \begin{equation}\label{eq2:a4spezial}
    a_4 = 2 a_1 (a_1+a_2)^{-1} a_2
\end{equation}
provided that $a_3=0$ and $a_1+a_2\ne 0$. Likewise, the remaining exposition is
extremely brief. The reader merely is invited to repeat the proof from
\cite[pp.~192--194]{schrei+s-35a} in order to verify that the given harmonicity
preserver of $L$ gives rise to a bijection $\sigma\colon K\to K$ satisfying
equations \eqref{eq:I} and \eqref{eq:II}.
\par
We present here a more detailed description, because we want to illustrate how
the definition of a semi-automorphism arises by following Ancochea's advise.
First, a proof of \eqref{eq2:a4spezial} can be done using
\begin{equation*}
    a_4 = \left( a_1^{-1} + a_2^{-1} \right)^{-1}
            \left( a_1^{-1}a_1 + a_2^{-1}a_2 \right)
        = 2 \left( a_2^{-1} (a_2+a_1) a_1^{-1} \right)^{-1} .
\end{equation*}
Next, in order to establish \eqref{eq:I}, it is enough to start on page 192 of
\cite{schrei+s-35a} before formula (15) and to repeat the reasoning from there
up to equation (20). Thereby, one has to reinterpret the cross ratios appearing
in \cite{schrei+s-35a}. Any cross ratio of the form $\cross(\infty,0,1,p)$ with
$p$ in $L\setminus\{\infty\}$ is to be understood as the corresponding element
$x_2^{-1}x_1\in K$ according to \eqref{eq:identif}; cross ratios equal to
$-1\in K$ characterise harmonic quadruples. The transfer from
\cite{schrei+s-35a} of the proof for \eqref{eq:II} appeared at the first sight
impossible to the author, since the commutativity of the ground field is used
at the very beginning of the reasoning in the middle of page 193. Fortunately,
it suffices to start reading on that page a few lines further down at equation
(21) taking into account that the fraction
\begin{equation*}
    \frac{2ab}{a+b}
\end{equation*}
appearing there has to be replaced by the right hand side of our equation
\eqref{eq2:a4spezial}. In the same way several other fractions, which do not
make sense over a skew field, have to be rewritten appropriately. The next
formula that needs to be altered is
\begin{equation*}
    a  b = \frac12 \left((a+b)^2-a^2-b^2 \right),
\end{equation*}
since this requires commuting elements $a, b$. In our not necessarily
commutative field $K$ this formula has to be replaced with
\begin{equation*}
    a_1 a_2 + a_2 a_1 =  (a_1+a_2)^2-a_1^2-a_2^2
                        \mbox{~~for all~~} a_1,a_2\in K.
\end{equation*}
In this way one readily obtains the condition in \eqref{eq:II} rather than the
multiplicativity of $\sigma$ as in \cite{schrei+s-35a}.
\par
After the proof, Ancochea gave examples of harmonicity preservers. Remarkably
enough, he did not use formula \eqref{eq:mult} for this purpose. Instead, a
result of S.~Wachs \cite[p.~82]{wachs-36b} was quoted. It states that, for all
$a_1,a_2,a_3,a_4\in K$, the condition $\Harm(a_1,a_2,a_3,a_4)$ is equivalent to
\begin{equation}\label{eq2:harm}
    (a_2-a_4)^{-1}(a_2-a_3)(a_1-a_3)^{-1}(a_1-a_4) = -1 .
\end{equation}
By virtue of this characterisation it turned out that \emph{any automorphism
and any anti-automorphism\footnote{Ancochea's terminology reads \emph{direct
automorphism} and \emph{reciprocal automorphism}, respectively.} of $K$ extends
to a bijective harmonicity preserver of $L$ fixing the points $0$, $1$,
$\infty$}.
\par
In \cite[\S~2--3]{anco-42a}, Ancochea dealt with an arbitrary quaternion skew
field $Q$ of characteristic $\neq 2$. He showed that \emph{any
semi-automorphism of $Q$ is either an automorphism or an anti-automorphism of
$Q$.} We do not enter a detailed analysis here, since in his subsequent article
\cite{anco-47a} a more general result can be found. The geometric significance
is, of course, that \emph{any semi-automorphism of $Q$ gives rise to a
bijective harmonicity preserver.}
\par
The article \cite{anco-42a} closes with some remarks on the case when the
centre of a quaternion skew field $Q$, $\Char Q\neq 2$, admits only the
identical automorphism. Ancochea recalled that, by what nowadays is called the
\emph{Theorem of Skolem-Noether}, all automorphisms of $Q$ are inner
\cite{noet-33a} (or see \cite[Thm.~4.9]{jaco-89a}). Making use of
\cite[Teorema~fundamental]{anco-41a} or the analogous result in \cite{reid-30a}
he derived that over such a quaternion skew field \emph{any bijective
harmonicity preserver either is a projectivity or---loosely speaking---the
product of some projectivity and the mapping that takes any quaternion in $Q$
to its conjugate}. He also referred to the description of all \emph{continuous}
bijective harmonicity preservers of the projective line over the real
quaternions given by S.~Wachs \cite[pp.~108--109]{wachs-36b}; he thereby
stressed that his reasoning was based exclusively on algebraic tools. Finally,
he emphasised the difference with the situation over the field of complex
numbers. Indeed, in order to obtain there a ``similar'' result, one has to
assume continuity due to the existence of discontinuous automorphisms of the
complex number field \cite{kamk-27a} (or see \cite{kest-51a}, \cite{mori-35a},
\cite{segre-47a}).
\par
Summing up, Ancochea generalised the first part of von Staudt's theorem to the
general case of an arbitrary field $K$, commutative or not. The second part,
namely the question whether or not all semi-automorphism of $K$ arise in this
way, remained as an open problem in \cite{anco-42a}.

\subsection{On semi-automorphisms of division algebras}\label{subse:semiautos}
The article \cite{anco-47a}, which is the third and last in this series,
commences with a short summary of the basic notions and results from
\cite{anco-42a}. In doing so, a \emph{semi-automorphism} of a ring $R$ is
defined as a bijection $\sigma\colon R\to R$ possessing the properties
\eqref{eq:I} and \eqref{eq:II} up to a change of notation. Further below,
\emph{semi-isomorphisms} of rings are introduced in the same fashion. Like
Ancochea, we tacitly assume from now on any \emph{algebra} to be associative
and finite-dimensional over its centre (unless explicitly stated otherwise).
So, \emph{division algebra} means a field with finite dimension over its
centre. The main results read as follows:

\begin{thm*}[{\cite[Principal Theorem]{anco-47a}}]
Let $D$ be a division algebra of characteristic $\neq 2$. Then any
semi-automorphism $\sigma$ of $D$ is an automorphism or an anti-automorphism of
$D$.
\end{thm*}

\begin{thm*}[{\cite[von Staudt's Theorem]{anco-47a}}]
In a projective space over a division algebra of characteristic $\neq 2$, any
bijective harmonicity preserver of a line $L$ onto itself with three distinct
fixed points $0$, $1$, $\infty$ restricts to a semi-automorphism of the
division algebra $D:=L\setminus\{\infty\}$. Conversely, any semi-automorphism
of $D$ extends to a bijective harmonicity preserver of $L$ that fixes the
points of reference.
\end{thm*}
Ancochea based the proof of his Principal Theorem upon several auxiliary
results, which are stated below.

\begin{thm*}[{\cite[Theorem~1]{anco-47a}}]
Under any semi-automorphisms $\sigma$ of a division algebra $D$ (of arbitrary
characteristic) the centre $F$ of $D$ is invariant.
\end{thm*}
In order to show his Theorem~1, Ancochea followed the lines of \cite{anco-42a}
by concluding that for all $\gamma\in F$ the image $\gamma^\sigma$ commutes
with all elements of the form
\begin{equation}\label{eq:ab-ba}
    ab - ba \mbox{~~with arbitrary~~} a,b\in D.
\end{equation}
Thus, it was enough to establish:

\begin{lem*}[{\cite[Lemma]{anco-47a}}]
If $c\in D$ commutes with all elements $ab-ba$ as in \eqref{eq:ab-ba}, then $c$
belongs to the centre $F$ of $D$.
\end{lem*}
In \cite{anco-47a}, there are two different proofs of this Lemma. The first one
is subject to the extra assumption $\Char D\neq 2$, the second one applies for
arbitrary characteristic of $D$. As main tool in both proofs, Ancochea
considered for any fixed element $x\in D\setminus F$ the mapping sending a
variable element $a\in D$ to $ax -xa$. This mapping is an $F$-endomorphism of
the $F$-vector space $D$ and its image is therefore an $F$-subspace $M_x$ of
$D$. Taking into account the dimension of $M_x$ and using results about the
structure of division algebras from the books of A.~A.~Albert \cite{albe-39a}
(reprint \cite{albe-61a}) and B.~van~der~Waerden \cite{vand-40a} (various
reprints and translations), Ancochea verified that for any element $c\in
D\setminus F$ there is at least one non-commuting element of the form
\eqref{eq:ab-ba}.
\par
Returning to the proof of the Principal Theorem (taking into account that
$\Char D\neq 2$), Ancochea first observed that the given semi-automorphism
$\sigma$ restricts to an automorphism of the commutative field $F$.
Furthermore, he obtained $(\gamma a)^\sigma = \gamma^\sigma a^\sigma$ for all
$\gamma\in F$, $a\in D$, i.~e., $\sigma$ is a semilinear mapping of the
$F$-vector space $D$. Next, an auxiliary $F$-algebra $D'$ was defined and an
$F$-semilinear isomorphism $\iota\colon D\to D'$ was explicitly established. He
noted that $\sigma_1:=\sigma\iota$ is a semi-isomorphism of $D$ onto $D'$.
Also, due to the particular choice of $\iota$, this $\sigma_1$ is
\emph{$F$-linear}.\footnote{In \cite{anco-47a} such a mapping is called a
semi-isomorphism \emph{over $F$}.} Ancochea's conclusion says that
\emph{$\sigma$ will be an (anti-)automorphism of $D$ precisely when $\sigma_1$
is an (anti-)isomorphism $D\to D'$}.
\par
After these technical preliminaries, Ancochea considered a commutative field
$K$ that is a finite algebraic extension of $F$ and a splitting field of $D$
and $D'$. Then he extended $\sigma_1$ to a $K$-linear semi-isomorphism
$\sigma'\colon D_K \to D_K'$, where $D_K:=K\otimes_F D$ and $D_K':=K\otimes_F
D'$. From $D_K$ and $D_K'$ being full matrix algebras of the same
$K$-dimension, followed their being $K$-isomorphic. Consequently, Ancochea
continued by assuming $\sigma_1$ to be a $K$-linear semi-automorphism of the
algebra $(\cM_n)_K$ of all $n\times n$ matrices over $K$ (for some positive
integer $n$). \emph{Evidently, $\sigma$ will be an (anti-)automorphism of $D$
if $\sigma'$ is an (anti-)automorphism.} In this way the proof of the Principal
Theorem was reduced to showing Theorem~2 below, which is important for its own
sake.

\begin{thm*}[{\cite[Theorem~2]{anco-47a}}]
Every $K$-linear semi-automorphism $\sigma$ of a full matrix algebra
$(\cM_n)_K$, where $K$ is a commutative field and $n$ is a positive integer, is
an automorphism or an anti-automorphism.
\end{thm*}
Ancochea's proof of Theorem~2 uses that $\sigma$ preserves orthogonal
idempotents of rank $1$. Let us write $e_{ij} \in (\cM_n)_K$ for the matrix
whose $(i,j)$-entry equals $1$, whereas all other entries vanish. Then it can
be assumed, without loss of generality, that $\sigma$ fixes the idempotents
$e_{ii}$, $i=1,2,\ldots, n$. The images of the matrices $e_{ij}$, where $i\neq
j$, turned out to satisfy $e_{ij}^\sigma = \alpha_{ij}e_{ij} +
\beta_{ij}e_{ji}$ with coefficients $\alpha_{ij},\beta_{ij}\in K$ subject to
$\alpha_{ij}\beta_{ij}=0$. This intermediate result led to two mutually
exclusive cases: either all $\beta_{ij}=0$ vanish or all $\alpha_{ij}$ vanish.
In the first case, a short calculation showed that $\sigma$ is an automorphism.
In the second case, $\sigma$ turned out to be an anti-automorphism, as
required.
\par
The version of von Staudt's Theorem in \cite{anco-47a} was then an immediate
consequence of \cite[Principal Theorem]{anco-47a} and the results from
\cite{anco-42a}.
\par
Ancochea closed his article with two important observations. On the one hand,
he sketched that the Principal Theorem remains true if $D$ is a \emph{simple
algebra of characteristic $\neq 2$}. He also claimed that the Principal Theorem
remained valid for semisimple algebras of characteristic $\neq 2$. However, the
last statement needs to be altered. We shall come across a correct version on
page~\pageref{page:correction}. On the other hand, he pointed out neat links to
the work of P.~Jordan, J.~von Neumann, and E.~Wigner \cite{jord-33a},
\cite{jord+v+w-34a}. Given an associative, but not commutative $F$-algebra $D$,
$\Char F\neq 2$, one obtains a commutative, but not associative algebra $D^+$
by maintaining the addition in $D$ and defining a new product $a\circ
b:=\frac12(ab+ba)$ for all $a,b\in D$. Obviously, \emph{the semi-automorphisms
of $D$ are precisely the automorphisms of $D^+$.} (Such an algebra $D^+$ has
been given the name \emph{special Jordan algebra} \cite[p.~178]{braun+k-66a}.)

\section{Rounding off Ancochea's results (1947--1953)}\label{se:rounding}
In 1947, I.~Kaplansky \cite{kapl-47a} mentioned that Ancochea's theorem on the
semi-automorphisms of a simple algebra \cite{anco-47a} fails for characteristic
two, because condition \eqref{eq:II} looses most of its strength.\footnote{For
example, if $K$ is a commutative field of characteristic two and $\sigma\colon
K\to K$ is an additive automorphism of $K$, then \eqref{eq:II} is trivially
satisfied for all $x,y\in K$ due to $(xy)^\sigma+(yx)^\sigma=0=x^\sigma
y^\sigma + y^\sigma x^\sigma$. If, furthermore, $K$ has more than two elements,
then there is a choice of $\sigma$ that does not fix $1\in K$ and therefore
cannot be a multiplicative automorphism.} In order to overcome this phenomenon,
Kaplansky modified the definition of a semi-isomorphism. We stick here to the
slightly different version that was given three years later by N.~I.~Jacobson
and C.~E.~Rickart \cite{jaco+r-50a}. Given (associative) rings $R$ and $R'$ a
mapping $\sigma\colon R\to R'$ is called a \emph{Jordan homomorphism} or, in
other words, a \emph{semi-homomorphism}\footnote{Both names are currently used
in the literature. We decided to switch freely between these names at our own
discretion.} if it satisfies for all $x,y\in R$ the following three conditions:
\begin{align}
    (x+y)^\sigma &= x^\sigma + y^\sigma,\label{eq:Ia} \\
    (x^2)^\sigma &= (x^\sigma)^2 ,\label{eq:IIa} \\
    (xyx)^\sigma &= x^\sigma y^\sigma x^\sigma .\label{eq:IIb}
\end{align}
Conditions \eqref{eq:I} and \eqref{eq:Ia} are identical. Any semi-homomorphism
$\sigma\colon R\to R'$ (in the above sense) satisfies Ancochea's condition
\eqref{eq:II} for all $x,y\in R$. This is immediately seen, using
\eqref{eq:IIa}, from $ xy+yx = (x+y)^2 - x^2 -y^2$. Conversely, any mapping
$\sigma\colon R\to R'$ satisfying Ancochea's conditions \eqref{eq:I} and
\eqref{eq:II} for all $x,y\in R$ is a semi-homomorphism provided that $R'$ is a
\emph{$2$-torsion free ring}, i.~e., $R'$ has no elements of additive order
two; see \cite[Lemma~1]{kapl-47a}, where Kaplansky used the identity
\begin{equation*}
    2xyx = 4(x+y)^3 - (x+2y)^3 - 3x^3 +4y^3 -2(x^2y+yx^2) .
\end{equation*}
\par
When dealing with \emph{unital rings} $R$ and $R'$ it is common to replace
\eqref{eq:IIa} with
\begin{equation}\label{eq:IIc}
    1^\sigma = 1'
\end{equation}
in the definition of a semi-homomorphism. As a consequence, with $y:=1$ in
\eqref{eq:IIb}, one obtains \eqref{eq:IIa} for all $x\in R$.
\par
Based on these observations and some lemmas, \cite{kapl-47a} comprises three
main results. In Theorem~1, any semi-isomorphism $\sigma\colon R\to R'$ is
considered, where $R$ and $R'$ are semisimple rings with unity. It is shown
that \emph{the restriction of $\sigma$ to the centre of $R$ yields an
isomorphism onto the centre of $R'$}. Theorem~2 describes the semi-isomorphisms
between simple algebras. \emph{Any mapping of this form turns out to be either
an isomorphism or an anti-isomorphism.} In Theorem~3, the erroneous statement
from \cite{anco-47a} about semi-automorphisms of semisimple algebras is
corrected and extended to a wider class of rings.
\label{page:correction}\emph{Let any semi-isomorphism $\sigma\colon R\to R'$ be
given, where $R$ and $R'$ are direct sums of simple algebras. Then the simple
components of $R$ and $R'$ may be paired off in such a way that the given
mapping is an isomorphism or an anti-isomorphism of each
pair.}\footnote{\label{fn:sum}Theorem~3 provides an easy way of constructing
semi-homomorphisms that are neither homomorphisms nor anti-homomorphisms. Take,
for example, the ring $\cM_n(F)$ of $n\times n$ matrices, $n\geq 2$, over any
commutative field $F$ and consider the mapping of the ring $\cM_n(F)\oplus
\cM_n(F)$ onto itself that sends any matrix pair $(A_1,A_2)$ to
$(A_1,A_2^{\top})$, where $A_2^{\top}$ denotes the transpose of $A_2$.} The
proofs of these theorems follow in part Ancochea's approach from
\cite{anco-47a}.
\par
Next, the article of L.-K.~Hua \cite{hua-49a} closed gaps that were left open
by Ancochea. Hua adopted the above definition of a Jordan homomorphism for
unital rings. His Theorem~1 states that \emph{any Jordan homomorphism of a
field $K$ onto itself is either an automorphism or an anti-automorphism of
$K$.} The way of proof differs considerably from Ancochea's, as it depends on a
series of subtle algebraic manipulations involving inverses rather than any
structure theory. Theorem~2 in \cite{hua-49a} provides an improvement and the
missing converse of Ancochea's theorem from \cite{anco-42a}. \emph{Over a field
$K$ of characteristic $\neq 2$, any bijective harmonicity preserver of a
projective line $L$ onto itself with three distinct fixed points determines an
automorphism or an anti-automorphism of $K$. Conversely, any surjective Jordan
endomorphism of $K$ gives rise to a bijective harmonicity preserver of the line
$L$.} A detailed exposition of the conclusions from \cite{hua-49a} was
published by Hua in \cite{hua-52a}. Furthermore, he sketched there ties to the
so-called \emph{geometry of matrices}, where it is also possible to
characterise harmonicity preserves by entirely different methods.\footnote{We
refer to \cite{wan-96a} and \cite{huanglp-06a} for the further development in
this area.}
\par
Let us return to N.~I.~Jacobson and C.~E.~Rickart, whose article
\cite{jaco+r-50a} contains the following theorem. \emph{Any Jordan homomorphism
of an arbitrary ring into an integral domain is either a homomorphism or an
anti-homomorphism}. Other results in \cite{jaco+r-50a} provide necessary
conditions for a Jordan homomorphism $\sigma\colon R\to R'$ to be the \emph{sum
of a homomorphism and an anti-homomorphism}. This is to mean that there are
ideals $R_1'$ and $R_2'$ of $R'$ with $R'=R_1'\oplus R_2'$, a homomorphism
$\sigma_1\colon R\to R_1'$ and an anti-homomorphism $\sigma_2\colon R\to R_2'$
such that $\sigma=\sigma_1+\sigma_2$. (The example in footnote \ref{fn:sum} is
of this kind.) This entails the existence of Jordan homomorphisms that are
neither a homomorphism nor an anti-homomorphism.
\par
The short communications \cite{herst-53a} and \cite{herst-53b} authored by
I.~N.~Herstein contain an alternative proof of Ancochea's Lemma from
\cite{anco-47a} on the elements of a division algebra $D$ commuting with all
elements as in \eqref{eq:ab-ba}.

\section{The impact of Ancochea's work}

\subsection{Applications in algebra and geometry}
There is a wealth of articles deploying results of Ancochea (primarily those
from \cite{anco-47a}) in order to solve a variety of problems in algebra. The
following list, which is in chronological order, comprises publications from
1959 up to 2005:
M.~Gerstenhaber \cite{gerst-59a}, 
N.~Jacobson \cite{jaco-59a}, 
H.-J.~Hoehnke \cite{hoeh-67a}, 
M.~Ra\"{\i}s \cite{rais-72a}, 
C.~V.~Devapakkiam \cite{deva+r-76a}, 
M.-A.~Knus \cite{knus+o-77a}, 
H.~F.~de~Groote \cite{degr-78a}, 
M.~Ra\"{\i}s \cite{rais-87a}, 
M.~O'Ryan and D.~B.~Shapiro \cite{oryan+s-96a}, 
R.~Parimala, R.~Sridharan, and M.~L.~Thakur \cite{pari+s+t-98a}, 
L.~Grunenfelder, T.~Ko\v{s}ir, M.~Omladi\v{c} and H.~Radjavi
\cite{grun+k+o+r-02a}, 
M.~A.~Chebotar, W.-F.~Ke and P.-H.~Lee \cite{cheb+k+l-05a}, 
M.~A.~Chebotar, W.-F.~Ke, P.-H.~Lee and L.-S.~Shiao \cite{cheb+k+l+s-05a}.
\par
Ancochea's article \cite{anco-42a} is one of the main sources in J.~Bilo's
monograph \cite{bilo-49b} about \emph{projective geometry over the real
quaternions}. Also, when dealing there with projectivities, a short comment
about the findings and methods used in \cite{anco-41a} is made on page~60. The
characterisation of harmonicity preservers (in Hua's version \cite{hua-49a}),
was used by W.~Benz \cite{benz-69a}, \cite[p.~175]{benz-73a},
\cite[p.~346]{benz-73a} in order to determine all isomorphisms of certain
\emph{M\"{o}bius geometries}. K.~List \cite{list-01a} made use of the same result
when characterising \emph{orthogonality preserving transformations} on the line
set of a three-dimensional \emph{hyperbolic space}.

\subsection{Semi-homomorphisms and their generalisations}
The outcomes of Ancochea together with the contributions by others (from
1947--1953) are at the beginning of a long series of articles. The results
proved there frequently read that---under certain extra conditions on the rings
$R$ and $R'$---\emph{any Jordan homomorphism $\sigma\colon R\to R'$ is either a
homomorphism or an anti-homomorphism.} By relaxing these extra conditions, the
conclusion often says that $\sigma$ is \emph{the sum of a homomorphism and an
anti-homomorphism}. We present a short overview and sketch various
generalisations.
\par
The result from \cite{jaco+r-50a}, saying that any Jordan homomorphism into an
integral domain is either a homomorphism or an anti-homomorphism, was shown
independently by L.-K.~Hua (under slightly stronger assumptions) in
\cite{hua-52a} and \cite{hua-53a}. In 1956, I.~N.~Herstein \cite{herst-56a}
extended these results to \emph{surjective} Jordan homomorphisms onto
\emph{prime rings}\footnote{A prime ring $R$ is one in which $aRb = 0$ implies
that $a = 0$ or $b = 0$.} with characteristic greater than three. K.~Yamaguti
\cite{yama-57a} established essentially the same theorem one year later.
M.~F.~Smiley \cite{smil-57a} gave a new proof covering also the case of
characteristic three and, by adopting Kaplansky's definition of a Jordan
homomorphism \cite{kapl-47a}, extended Herstein's theorem to the missing case
of characteristic two. The next steps were taken between 1979 and 1989 by
W.~E.~Baxter and W.~S.~Martindale \cite{baxt+m-79a}, M.~Bre\v{s}ar
\cite{bres-89a}, \cite{bres-91a}, who went over from prime to \emph{semiprime
rings} \cite[p.~158]{lam-01a}.
\par
Let us go back to the year 1948, when N.~Jacobson \cite{jaco-48a} took up a
remark of Ancochea by saying that the statements from \cite{anco-47a} can be
seen as results about the isomorphisms of the non-associative special Jordan
rings determined by the given rings. Based on this observation, Jacobson
initiated the study of isomorphisms between various \emph{Jordan subsystems} of
(associative) rings. In this way he gave also new proofs for some results from
\cite{anco-47a}. He also noted that, by replacing the plus sign with a minus in
Ancochea's condition \eqref{eq:II}, one arrives at another kind of
``semi-isomorphism'' of associative rings, which is related to the theory of
\emph{Lie rings}.\footnote{Several of our bibliographical items deal also with
this topic. Further information may be retrieved from \cite{beid+c-01a} and
\cite{beid+m+c-04a}.} His work led to a series of results stating that
\emph{any Jordan homomorphism between certain Jordan subsystems can be lifted
to a homomorphism or an anti-homomorphism of the ambient rings}. In 1949,
F.~D.~Jacobson and N.~Jacobson \cite{jaco+j-49a} determined all embeddings of a
special Jordan algebra in associative algebras over a commutative field of
characteristic $\neq 2$, thereby regaining a theorem from \cite{anco-47a}.
N.~Jacobson and C.~E.~Rickart \cite{jaco+r-52a} considered the Jordan subsystem
of symmetric elements of a ring with involution. Their work was later extended
by
W.~S.~Martindale \cite{mart-67a}, 
L.~A.~Lagutina \cite{lagu-88a}, 
K.~McCrimmon \cite{mccr-89a}, 
W.~S.~Martindale \cite{mart-90a}, 
K.~I.~Beidar and M.~A.~Chebotar \cite{beid+c-00a}. 
For summaries and extensive bibliographies we refer to the work
of 
I.~N.~Herstein \cite{herst-61a}, 
N.~Jacobson \cite{jaco-68a}, 
I.~N.~Herstein \cite{herst-69a}, 
N.~Jacobson \cite{jaco-76a}, 
W.~S.~Martindale \cite{mart-77a}, 
R.~P.~Sullivan \cite{sull-83a}, 
K.~McCrimmon \cite{mccr-89a}.
\par
A detailed analysis of the Jordan homomorphisms of the ring $\cT_r(R)$ of
\emph{upper triangular $r\times r$-matrices} with entries in a ring $R$
commenced in 1998 with the article \cite{moln+s-98a} of L.~Moln\'{a}r and
P.~\v{S}emrl on \emph{linear preserver problems}. They described all
\emph{linear} Jordan automorphism of $\cT_r(\bC)$, where $\bC$ denotes the
field of complex numbers. Shortly afterwards, K.~I.~Beidar, M.~Bre\v{s}ar,
M.~A.~Chebotar \cite{beid+b+c-00a} switched to triangular matrices over a
unital commutative ring $R$. One of their results is in the spirit of Ancochea
and illustrates the extra features arising from idempotent ring elements.
\emph{Let $R$ be a $2$-torsion-free commutative ring with identity. Then $R$
contains no idempotents except $0$ and $1$ if, and only if, every Jordan
isomorphism of $\cT_r(R)$, $r\geq 2$, onto an arbitrary algebra over $R$ is
either an isomorphism or an anti-isomorphism.} Further contributions (with
varying assumptions on the ring $R$) have been given between 2001 and 2014:
X.~M.~Tang, C.~G.~Cao and X.~Zhang \cite{tang+c+z-01a}, 
X.~T.~Wang and H.~You \cite{wang+y-04a}, 
D.~Benkovi\v{c} \cite{benk-05a}, 
T.~L.~Wong \cite{wong-05a}, 
C.-K.~Liu and W.-Y.~Tsai \cite{liu+t-07a}, 
X.~T.~Wang \cite{wang-07a}, 
H.~M.~Yao, C.~G.~Cao and X.~Zhang \cite{yao+c+z-09a}, 
L.-P.~Huang \cite{huanglp-10b}, 
X.~T.~Wang and Y.~M.~Li \cite{wang+l-10a}, 
Y.~Wang and Y.~Wang \cite{wang+w-13a} 
(erratum by Y.~Du, Y.~Wang and Y.~Wang \cite{du+w+w-14a}), 
Y.~Du and Y.~Wang \cite{du+w-14a}.
\par
Jordan homomorphisms between other kinds of rings have been thoroughly
investigated. We refer to
J.~H.~Zhang \cite{zhang-02a}, 
K.~I.~Beidar, M.~Bre\v{s}ar and M.~A.~Chebotar \cite{beid+b+c-03a}, 
F.~Lu \cite{lu-03a}, 
F.~Lu and T.-L.~Wong \cite{wong-05a}, 
A.~L.~Yang and J.~H.~Zhang \cite{yang+z-08a}, 
F.~Kuzucuo\u{g}lu and V.~M.~Levchuk \cite{kuzu+l-10a}, 
J.~Alaminos, J.~Extremera and A.~R.~Villena \cite{alam+e+v-12a}, 
E.~Akkurt, M.~Akkurt and G.~P.~Barker \cite{akku+a+b-15a}, 
A.~L.~Yang \cite{yang-15a}, 
R.~Brusamarello, E.~Z.~Fornaroli and M.~Khrypchenko \cite{brus+f+k-18a}.
\par
There are many articles dealing with mappings between rings that satisfy a
\emph{relaxed version} of the properties defining a Jordan homomorphism
$\sigma\colon R\to R'$ or some ``similar'' functional equations. Two types have
found particular interest.
First, \emph{$n$-Jordan mappings} are additive and satisfy $(x^n)^\sigma=
(x^\sigma)^n$ for all $x\in R$ and some fixed integer $n\geq 2$. These mappings
have been investigated by
I.~N.~Herstein \cite{herst-56a}, \cite{herst-61a}, \cite{herst-69a}, 
M.~Bre\v{s}ar, W.~S.~Martindale and C.~R.~Miers \cite{bres+m+m-98a}. 
Second, we mention \emph{Jordan triple homomorphisms}, which are characterised
by the conditions in \eqref{eq:Ia} and \eqref{eq:IIb} and are the topic of
articles by 
K.~Yamaguti \cite{yama-57a}, 
I.~N.~Herstein \cite{herst-67a}, 
M.~Bre\v{s}ar \cite{bres-89a}, \cite{bres-91a}, 
L.~A.~Lagutina \cite{lagu-91a}, 
F.~Lu \cite{lu-08a}.
Other contributions in this spirit came between 1953 and 1999 from 
J.~K.~Goldhaber \cite{gold-53a}, 
E.~Artin \cite[pp.~37--40]{artin-57a}, I.~N.~Herstein and E.~Kleinfeld
\cite{herst+k-60a}, 
M.~F.~Smiley \cite{smil-61a}, 
D.~W.~Barnes \cite{barn-66a}, 
K.~McCrimmon \cite{mccr-71a}, 
S.~A.~Huq \cite{huq-87a}, 
R.~Awtar \cite{awtar-88a}, K.~I.~Beidar and Y.~Fong \cite{beid+f-99a}, 
and after the year 2000 from
K.~I.~Beidar, S.-C.~Chang, M.~A.~Chebotar and Y.~Fong \cite{beid+c+c+f-00a}, 
M.~Bre\v{s}ar \cite{bres-00a}, 
M.~A.~Chebotar, F.-W.~Ke, P.-H.~Lee and L.-S.~Shiao \cite{cheb+k+l+s-05a}, 
M.~Bre{\v{s}}ar, M.~A.~Chebotar and W.~S.~Martindale \cite{bres+c+m-07a}, 
A.~K.~Faraj, A.~H.~Majeed, C.~Haetinger and N.~u.~Rehman
\cite{faraj+m+h+n-14a}. 
\par
Yet another way of generalising the notion of semi-homomorphism consists in
considering algebraic structures other than rings. Already in 1951, F.~Dinkines
\cite{dink-51a} 
introduced \emph{semi-homomorphisms of groups} by following condition
\eqref{eq:IIb}. Her results were extended by
I.~N.~Herstein and M.~F.~Ruchte \cite{herst+r-58a}, 
I.~N.~Herstein \cite{herst-68a}, 
K.~I.~Beidar, Y.~Fong, W.-F.~Ke, W.-R.~Wu \cite{beid+f+k+w-99a}. 
In the 1980s, R.~P.~Sullivan \cite{sull-83a}, \cite{sull-85a} considered also
\emph{semigroups} and introduced even more general \emph{half-automorphisms}.
Other generalisations deal with mappings between \emph{alternative division
rings}, \emph{near rings} and \emph{semirings}; see
M.~F.~Smiley \cite{smil-62a}, 
K.~C.~Smith and L.~van~Wyk \cite{smith+v-92a}, 
K.~I.~Beidar, Y.~Fong, W.-F.~Ke and W.-R.~Wu \cite{beid+f+k+w-99a}, 
S.~Shafiq and M.~Aslam \cite{shaf+a-17a}, 
B.~L.~M.~Ferreira and R.~N.~Ferreira \cite{ferr+f-19a}.

\subsection{Harmonicity preservers}
The complete description of all bijective harmonicity preservers of the
projective line $L$ over a field $K$ with characteristic $\neq 2$ has become a
standard topic in textbooks on projective geometry. We confine ourselves to
quoting several books from the 1950s.\footnote{Further references and
historical remarks are contained in \cite[pp.~57--58]{karz+k-88a}.} The first
exposition appears to be the one of R.~Baer \cite[p.~78]{baer-52a}. He
considered a more general setting, namely bijections of $L$ that preserve an
arbitrary cross ratio $d\neq 0,1$ lying in the centre of $K$ without any
restriction on the characteristic of $K$. Again, it was enough to treat the
case of such a bijection with three distinct fixed points $0$, $1$, $\infty$.
Baer showed that such a bijection gives rise to an automorphism or an
anti-automorphism of $K=L\setminus\{\infty\}$ fixing the element $d\in K$ and
vice versa. He thereby generalised also an outcome of A.~J.~Hoffman
\cite{hoff-51a}, who had obtained the same kind of result for a commutative
field $K$. An alternative proof of Baer's theorem can be found in the book of
G.~Pickert \cite[p.~121]{pick-55a} from 1955 (second edition \cite{pick-75a}).
E.~Artin \cite[pp.~84--85]{artin-57a} reproved the original result about
harmonicity preservers in an elegant way using an alternative characterisation
of semi-automorphisms \cite[pp.~37--40]{artin-57a}.
\par
The uniqueness of the fourth harmonic point (see page
\pageref{page:uniqueness}) holds more generally in projective \emph{Moufang
planes}, which satisfy a weaker form of Desargues' Theorem \cite{mend-56a},
\cite{pick-55a}, \cite{pick-75a}. These projective planes admit coordinates
from an \emph{alternative division ring} and the bijective harmonicity
preservers between their lines admit an algebraic description in terms of
Jordan isomorphisms. This investigation was initiated in the 1950s by V.~Havel
\cite{havel-55a}, \cite{havel-55b}, \cite{havel-57a} and continued by
M.~F.~Smiley \cite{smil-62a}, J.~van~Buggenhout \cite{vanb-69a}, J.~C.~Ferrar
\cite{ferr-81a}. Generalisations (in terms of other cross ratios) were given by
A.~Schleiermacher \cite{schl-65a}, H.~Schaeffer \cite{schae-81a}, A.~Blunck
\cite{blunck-91a}.
\par
In \cite{havel-68a}, V.~Havel considered harmonicity preservers in certain
\emph{translation planes}. The work of W.~Bertram \cite{bert-03a} on
harmonicity preserving mappings of some other geometric structures also
deserves mention.
\par
An area of ongoing research is the study of harmonicity preservers between
\emph{projective lines over unital rings} (with $1\neq 0$). Given such a ring
$R$ one considers any free left $R$-module $M$ that has at least one basis with
two elements. The point set of the projective line $\bP(M)$ is the set of all
cyclic submodules $Ru$ of $M$ for which there is a $v\in M$ such that $(u,v)$
is a basis of $M$ \cite{blunck+he-05a}, \cite{herz-95a}. Points $p_1,p_2$ of
$\bP(M)$ are called \emph{distant} (or \emph{non-neighbouring}) if there exists
a basis $(u_1,u_2)$ of $M$ with $p_1=Ru_1$ and $p_2=Ru_2$. This distant
relation is symmetric and anti-reflexive, and it turns the point set of
$\bP(M)$ into the so-called \emph{distant graph} of
$\bP(M)$.\footnote{Precisely when $R$ is a field, ``being distant'' means the
same as ``being distinct''.}
\par
Notions like \emph{harmonic quadruple} and \emph{cross ratio} are defined on
$\bP(M)$ in almost same way as on $\bP(V)$ when $V$ is a left $K$-vector space.
But, instead of ``distinct points'' of $\bP(V)$ one has to consider ``distant
points'' of $\bP(M)$. Furthermore, upon choosing a basis $(w_1,w_2)$ of $M$,
one obtains $0:=Rw_2$, $1:=R(w_1+w_2)$ and $\infty:=Rw_1$ as three \emph{points
of reference} on $\bP(M)$. It is straightforward then to identify the points of
$\bP(M)$ that are distant to $\infty$ with the elements of $R$ as Ancochea did
in \cite{anco-42a}. For doing so, it is enough to replace in \eqref{eq:identif}
the field $K$ by the ring $R$ and to assume that $x_2\in R$ is invertible.
However, unless $R$ is a field, the ``rest'' of the projective line contains
apart from $\infty$ many other ``points at infinity''.
\par
There is a widespread literature on \emph{harmonicity preservers}, which are
defined in the same manner as on page \pageref{page:h-preserver}, under varying
assumptions on the underlying unital rings $R$ and $R'$. Below we collect the
relevant contributions. A common feature in all of them is that \emph{$2=1+1$
has to be invertible in $R$}.
\par
Let a harmonicity preserver $\lambda\colon \bP(M)\to\bP(M')$ be given. Then,
after the identification of a subset of $\bP(M)$ with $R$ and an analogous
identification in $\bP(M')$, \emph{$\lambda$ restricts to a Jordan homomorphism
$R\to R'$ provided that $R$ contains ``sufficiently many'' invertible
elements.} A proof of this result can be derived from 
B.~V.~Limaye and N.~B.~Limaye \cite{lima+l-77b}, 
despite the fact that their work from 1977 is mainly about \emph{commutative}
rings. A formal proof under slightly weaker assumptions was given by the 
author \cite{havl-15b}. 
Already in 1971, 
N.~B. Limaye \cite{lima-71a}, \cite{lima-72a} 
proved a version for \emph{commutative local} rings \cite[pp.~280f.]{lam-01a}
and then for \emph{commutative semilocal} rings \cite[p.~296]{lam-01a}; 
H.~Schaeffer \cite{schae-74a}, 
B.~V.~Limaye and N.~B.~Limaye
\cite{lima+l-77b}, 
B.~R.~McDonald \cite{mcdon-81a} 
treated also the case of \emph{commutative} rings.
\par
The converse problem, like before, is to decide whether or not \emph{any Jordan
homomorphism $\sigma\colon R\to R'$ gives rise to a harmonicity preserver}.
After choosing bases $w_1, w_2$ of $M$ and $w_1',w_2'$ of $M'$ it is tempting
to define a mapping
\begin{equation}\label{eq:nosolution}
    M \to M' \colon x_1w_1 + x_2w_2 \mapsto  x_1^\sigma w_1' + x_2^\sigma w_2'
    \mbox{~~with~~} x_1,x_2\in R .
\end{equation}
However, this mapping will in general not give rise to a mapping
$\bP(M)\to\bP(M')$, let alone its being a harmonicity preserver. Nevertheless,
for commutative rings the approach in \eqref{eq:nosolution} does work, since a
Jordan homomorphism of $R$ in $R'$ is nothing but a homomorphism. For \emph{not
necessarily commutative} rings the situation is much more involved, due to the
possibly large number of ``points at infinity''. B.~V.~Limaye and N.~B.~Limaye
\cite{lima+l-77a} (erratum \cite{lima+l-77c}) gave an affirmative answer to the
problem for \emph{local} rings. They actually determined all bijections of the
projective line $\bP(M)$ onto itself such that all quadruples with a given
cross ratio $d$ go over to quadruples with a given cross ratio $d'$, where $d,
d'$ are elements in the centre of $R$ other that $0,1$. A.~Herzer
\cite{herz-87a} showed how to obtain well-defined point-to-point mappings from
certain Jordan homomorphisms. A breakthrough is due to C.~Bartolone
\cite{bart-89a}, who defined the desired mapping $\bP(M)\to\bP(M')$ under the
extra condition that $R$ is a ring of \emph{stable rank two}
\cite[p.~24]{blunck+he-05a}. In terms of the bases used in
\eqref{eq:nosolution}, his solution\footnote{The authors of \cite{herz-87a} and
\cite{bart-89a} had different aims and did not exhibit the preservation of
harmonicity in their publications.} reads
\begin{equation}\label{eq:bartolone}
    R\bigl( x w_1 + (1+xy)w_2\bigr) \mapsto
    R'\bigl(x^\sigma w_1' + (1'+x^\sigma y^\sigma) w_2'\bigr)
    \mbox{~~with~~} x,y\in R.
\end{equation}
A.~Blunck and the author treated in \cite{blunck+h-03b} the general case taking
into account that the distant graph on $\bP(M)$ need not be connected. It
turned out that \emph{a Jordan homomorphism $\sigma$ determines a harmonicity
preserving mapping only on that connected component of the distant graph on
$\bP(M)$ which contains the chosen points of reference}. (The formulas used for
this purpose arise from \eqref{eq:bartolone} in a recursive way.) As a
consequence, one may select arbitrarily \emph{one Jordan homomorphism $R\to R'$
per component} in order to create a harmonicity preserver $\bP(M)\to\bP(M')$.
\par
The material from the last two paragraphs forms the foundation for the version
of von Staudt's Theorem in \cite[Thm.~1]{havl-15b}. When dropping the
assumption on the existence of ``sufficiently many units'' in $R$, this theorem
fails. A lucid counterexample in terms of the polynomial ring in one
indeterminate over the real numbers was given by C.~Bartolone and F.~Di~Franco
\cite{bart+f-79a} already in 1979. They therefore initiated the study of
mappings that preserve \emph{generalised harmonic quadruples} and succeeded in
describing all such mappings for commutative rings; see also M.~Kulkarni
\cite{kulk-80a}, C.~Bartolone and F.~Bartolozzi \cite{bart+b-85a},
L.~Cirlincione and M.~Enea \cite{cirl+e-90a}, A.~A.~Lashkhi \cite{lash-89a},
\cite{lash-90a}, \cite{lash-97a}, D.~Chkhatarashvili \cite{chkh-98a}. Closely
related is the work of F.~Buekenhout \cite{buek-65a}, St.~P.~Cojan
\cite{cojan-85a}, D.~G.~James \cite{james-82a}, B.~Klotzek \cite{klotz-88a},
who characterised mappings that satisfy a rather weak form of cross ratio
preservation between projective lines over fields. The algebraic background of
their work is a re-coordinatisation of the domain projective line in terms of a
\emph{valuation ring} \cite[p.~1]{math-86a}.
\par
All things considered, Ancochea's semi-homomorphisms keep going strong. They
constitute the indispensable algebraic tool for describing harmonicity
preservers between projective lines over unital rings.


\footnotesize

\newcommand{\Dbar}{\makebox[0cm][c]{\hspace{-2.5ex}\raisebox{0.25ex}{-}}}\newcommand{\cprime}{$'$}

\noindent
Hans Havlicek\\
Institut f\"{u}r Diskrete Mathematik und Geometrie\\
Technische Universit\"{a}t\\
Wiedner Hauptstra{\ss}e 8--10/104\\
A-1040 Wien\\
Austria\\
\texttt{havlicek@geometrie.tuwien.ac.at}

\end{document}